\input amstex
\input xy
\xyoption{all}
\documentstyle{amsppt}
\document
\magnification=1200
\NoBlackBoxes
\nologo
\hoffset1.5cm
\voffset2cm
\pageheight {16cm}
\def\C{\Cal{C}}

\def\U{\Cal{U}}
\def\N{\Cal{N}}

\bigskip

\centerline{\bf NEURAL CODES AND HOMOTOPY TYPES: }

\medskip

\centerline{\bf MATHEMATICAL MODELS OF PLACE FIELD RECOGNITION}

\bigskip

\centerline{\bf Yuri I.~Manin}

\medskip

\centerline{\it Max--Planck--Institut f\"ur Mathematik, Bonn, Germany}
\bigskip

\hfill{ To Misha and Serezha, cordially}
\bigskip

{\it ABSTRACT. This note is a brief survey of some results of the recent collaboration of
neurobiologists and mathematicians dedicated to stimulus
reconstruction from neuronal spiking activity.
This collaboration, in particular, led to the consideration of binary codes
used by brain for encoding a stimuli domain such as a rodent's territory through the combinatorics
of its covering by local neighborhoods.

The survey is addressed to mathematicians (cf. [DeSch01]) and focusses on the idea that
stimuli spaces are represented by the relevant neural codes
as simplicial sets and thus encode say, the homotopy type of space
if local neighbourhoods are convex (see [CuIt08], [CuItVCYo13], [Yo14], [SiGh07]).}

\bigskip

{\bf 1. Introduction and summary.}  In this note, {\it a code $\Cal{C}$} in general means a   {\it set  of   words
in a finite alphabet.}  Finite sequences of code words encode  certain data,
and mathematical problems related to codes depend on the class of situations 
in which encoding/decoding of the relevant data occurs. For the restricted purposes
of this introduction, we will divide  these situations and the relevant problems
into three groups.

\smallskip

(i) {\it Codes and encoding/decoding procedures in cryptography} must efficiently
work in those situations where there is a challenge of unauthorised 
access to and/or falsification  of the encoded information, cf. [Ya00].

\smallskip

In the relevant group of mathematical problems, encoding/decoding programs
must be computationally easy, but {\it the choice of the actually used program}
by a pair (or a larger group) of users must itself be encoded in such a way,
that guessing it becomes computationally unfeasible for an unauthorised person/computer intervention.

\smallskip 

Since the notion of computational feasibility itself became formalised 
in a highly non--trivial way, in these studies the P/NP problem is often invoked.

\smallskip

(ii) {\it Error--correcting codes} must be constructed in such a way that
when the information encoded by a sequence of code words  is transmitted through a noisy channel, 
decoding algorithms could reconstruct the initial
text from its corrupted version, at least with high pro\-bability.
M.~Tsfasman and S.~Vladut were among the pioneers of
the creation and study of efficient algebraic geometric error--correcting codes
and coauthored an influential monograph about them that exists now
in several editions and translations, cf. [Vla\-NoTsfa07].

\smallskip

Interesting unsolved mathematical problems in the theory of
error--correcting codes are stated in terms of properties
of the set of all (unstructured or structured) codes of arbitrary length with a fixed
cardinality of the alphabet.
The study of mutually contradictory requirements of high transmission speed,
high probability of error--free correction, and computational feasibility
of encoding/decoding procedures, led to the notions
of {\it asymptotic bound for codes} and {\it isolated codes} (lying above the asymptotic bound):
see  [Vla\-NoTsfa07].
A recent progress is related with the understanding that if the Kolmogorov
complexity of the code is considered as its {\it energy},
the standard models of statistical physics furnish
an interpretation of the asymptotic bound as a phase transition curve:
see [Ma11] and [MaMar12].

\smallskip

(iii) Finally, {\it neural codes} to which this short survey is dedicated
are mathematical models suggested for the understanding
how brain copes with multiple tasks of orienting and navigating  in the world, cf. [CuIt08],
[CuItVCYo13], [Yo14], [GiIt14],  and references therein.

\smallskip

The starting point of discussion is the following remark: external stimuli,
say the spatial environment in which  an animal moves, is accessible for (or even created by)
a scientist. However, the brain of the animal must be able to reconstruct, say, a map
of this environment and the current position in it using only the action potentials (spikes) of the relevant cell groups.

\smallskip
The philosophy underlying construction of relevant mathematical models can be 
summarised as follows. In laboratory experiments designed to
study brain functions it is found that stimuli are naturally divided
into groups, and with each group a certain type of neural activity
is associated.

\smallskip

In [CuIt08] and [Yo14], it is postulated that a given type of stimuli can
be modelled via a topological, or metric
 {\it stimuli space $X$}. 
Furthermore,
brain reaction to a point in $X$ is modelled by spiking
activity of certain finite set of neurons $\N_X$. 

\smallskip

For example, whenever a rat in laboratory is moving in a restricted space,
with possible obstacles, $X$ is a  map of accessible territory,
whereas regions in which the animal finds himself, 
are marked by specific activities of {\it hippocampal place cells}: these cells 
constitute elements of $\N_X$.  

\smallskip
In an unrelated study   ([MaGaJoGoAshFraFri00]) of licensed London taxi drivers
 who had to pass an intense two--year
training for orientation in (the mental map of) London,
it was found that the total mass of hippocampal place cells
is considerably larger that in the control group.
(The authors worked with data from the pre--GPS tracker's  era).
This study is a strong argument for the assumption
that hippocampus carries orientation/navigation neural networks
in humans as well as in rodents.

\smallskip

Simplifying  and averaging observable spiking data obtained from neurons
in $\N_X$ at various points in $X$, the researchers ([CuIt08], [CuItVCYo13], [Yo14]) postulate the
following structures.

\smallskip
(a) The stimuli space $X$ is endowed with a certain covering $\{U_i\}$ by subsets called
{\it receptive fields.}

\smallskip

(b) Possible patterns of spiking activities  correspond to words 
belonging to the respective binary neural code. Each word in such code
is a map $w:\,\N_X\to \{0,1\}$ or equivalently, the
subset  $w^{-1}(1)\subset \N_X$ consisting of neurons that fire
at certain points of $X$ whereas   the remaining
neurons are inactive. 

\smallskip

(c) The correspondence between certain subsets of positions
in the stimuli space  and 
patterns of neural activity is given by the formula (1) below.

\smallskip

(d) During the period of time when an animal studies his territory,
(an appropriate group of cells in) his brain generates the code reflecting a map of it.
When this study evolves to the stage of active life in this territory,
sequences of code words encode paths in the territory.

\smallskip

The main problem which is addressed in [CuIt08] and [Yo14] is this:
{\it what properties of $X$ can be reconstructed if one knows only all occurring 
patterns of neuronal activities, that is mathematically, the binary code
$\C:= \C_X \subset  \{0,1\}^n$?} Here and below we fix an $X$ and number all neurons in $\N_X$
by $\{1,\dots ,n\}$ in order to simplify notation.

\smallskip

It turns out, in particular, that under some additional restrictions upon
$X$ and $\{U_i\}$, {\it the total homotopy type}
of $X$ can be inferred from  the relevant code $\C$.
In particular, finite sequences of words in $\C$ 
may be used  in order to model elements of the
{\it fundamental groupoid} of $X$: accessible paths from one neighbourhood in $X$
to another one. (Finding such paths of minimal length
used to be a constant challenge for brains of taxi drivers.)

\smallskip

The suggestion that brain uses and somehow encodes the combinatorics
of coverings (and more generally, diagrams in (poly)categories) was made already in 
[BrPo03].

\smallskip

{\it Plan of the paper.} In the section 2, I explain basic geometric structures
relating binary codes to topological objects that model stimuli spaces,
and algebraic--geometric structures useful for description of codes themselves.
The section 3 contains statements of results from [Yo14] showing how neural codes
describe homotopy types of stimuli spaces. Finally in section 4, I discuss
possible versions of this approach for psychological and semiotic studies of languages.

\medskip

{\bf 2.  Binary codes and related geometric objects.} Let now $\C$ be a binary code
as above. Define the {\it support
of a code word} $w=(\varepsilon_1,\dots ,\varepsilon_n)$ as
$\roman{supp}\,(w):= \{i\,|\varepsilon_i=1\} \subset \{1,\dots ,n\}$.
The correspondence between $\C$ and set of all supports
of words in $\C$ is a bijection, and we will sometimes identify 
$\C$ with the relevant subset of subsets in $\{1,\dots ,n\}$.

\smallskip

The binary code $\C$ is called {\it simplicial} one if for each $w\in \C$
and each $v\in \{0,1\}^n$ with $\roman{supp}\,(v) \subset \roman{supp}\,(w)$,
we have $v\in \C$. Clearly, arbitrary binary code $\C$ is contained in the unique
minimal simplicial code which we will denote $\overline{\C}$.

\smallskip

Following [Yo14], we will consider three types of geometric objects related to binary codes.

\smallskip

(i) {\it Simplicial set of a code.} The set of (supports) of all  words of a simplicial  code $\C$
has the natural structure of  a {\it simplicial set} (see e.~g.~ [GeMa03], sec. I.2).
To be more precise, this set of code words defines a {\it triangulated space }([GeMa03], sec I.1)
but the difference between the two notions is inessential in our context.
\smallskip

We denote the respective
topological space  $\Delta(\C )$.
\smallskip

If the code $\C$ is not simplicial, we put by definition $\Delta(\C ):=\Delta(\overline{\C} )$.
\smallskip

(ii) {\it Code of a  covering.} Let $X$ be a set endowed
with a finite family of  subsets $\U := \{ U_i\,|\,i=1,\dots ,n\}$. Define the binary code $\C(\U )$ (of length $n$)
of this family by the following condition ([Yo14], p.~9): 
$$
w\in \C(\U)\quad \roman{iff}\quad
\left(  \bigcap_{i\in \roman{supp}\,(w)}U_i\right) \setminus\left( \bigcup_{j\notin \roman{supp}\,(w)}U_j\right) \neq \emptyset.
\eqno (1)
 $$
We will usually assume that $X=\cup_{i=1}^n U_i.$
\smallskip

(iii) {\it Algebraic geometric description of a binary code.} Identify $\{0,1\}$ with $\bold{F}_2$: finite field with two
elements.  Consider $\{0,1\}^n$ as the set of  $\bold{F}_2$--points of coordinatized
$n$--dimensional affine space over $\bold{F}_2$. 
\smallskip
Any polynomial $P\in \bold{F}_2[x_1,\dots ,x_n]$ determines a code $\C_P$:
this code consists of all points of $\bold{F}_2^n$ at which $P$ vanishes.
\smallskip

Call $P$ a {\it reduced polynomial}, if each monomial with non--zero coefficient
in $P$ is a product of pairwise distinct variables $x_i$. 

\smallskip

An easy count (see e.~g.~[Ma99]) shows that the map $P\mapsto \C_P$ determines a bijection between reduced
polynomials and binary codes of length $n$. In other words, each code $\C$ consists of $\bold{F}_2$--points
of a unique hypersurface determined by a reduced polynomial over  $\bold{F}_2$.
The total ideal of polynomials vanishing at all points of $\C_P$ is generated by $P$ and
$x_i^2-x_i, i=1,\dots , n.$

\smallskip

Such and similar representations of a binary code are used in [Yo14] in order to devise
algorithms for extracting geometric information from neural codes related to stimuli spaces.

\medskip

{\bf 3. Encoding topology of receptive fields.} Given a covering $\U:\ X=\cup U_i$ as above,
we can construct its {\it nerve}:  the triangulated space $N(\U)$, whose vertices $u_i$ bijectively
correspond to $U_i$, whereas 
$$
k+1\ \roman{vertices}\ u_{i_0},\dots , u_{i_k}\  \roman{span\ a}\
k-\roman{simplex\ iff}\ \cap_{j=0}^k U_{i_j}\ne \emptyset.
\eqno(2)
$$

\smallskip

{\bf 3.1. Proposition.} ([Ha02], p. 459, Corollary 4G.3.) {\it Assume that $X$ is a paracompact 
topological space, $\U$ its open covering, and all non--empty finite
intersections of  $U_i$ are contractible. Then homotopy types of $X$ and $N(\U )$ are canonically
equivalent.}

\medskip

Notice that the contractibility assumption  on $U_i$'s in the orientation problems
models the fact that visual field at any point of space is covered by the direct lines of
vision in all directions.

\medskip

{\bf 3.2. Proposition.} ([Yo14], p.~13.) {\it With the same assumptions as above,
$N(\U )$ can be canonically identified with the triangulated space $\Delta (\C(\U )$.
\smallskip
Therefore, the total homotopy type of $X$ is encoded by $\C$.}

\smallskip

A check of  the statement $N(\U ) = \Delta (\C(\U ))$ is omitted in [Yo14],
but the author has kindly sent me an easy combinatorial argument. I will reproduce
it here in order to register that it {\it does not use any topological structure} of $X$:
we may simply assume that $X$ is a set represented as a union of its subsets $\U_i$
to which we apply definitions (1) and (2).

\smallskip

The inclusion $\Delta(\C(\U))\subset N(\U)$ follows from the fact that if
$\sigma\in \Delta(\C(\U))$, then there is $w\in \C(\U)$ whose support contains
$\sigma$, so that applying (1) to $w$ we get (2).

\smallskip

Conversely, $N(U)\subset  \Delta(\C(\U))$. In fact, consider $\sigma\in N(\U)$.
Choose $p\in \cap_{i\in \sigma}U_i$: it exists in view of (2). Put
$\tau = \{ j\in \{1,\dots , n\}\,|\, p\in U_j\}.$ Then $\sigma\subset \tau$
and applying (1) we see that $\sigma$ belongs to the simplicial
completion of $\C(\U)$ and hence $\sigma\in \Delta(\C(\U))$.

\smallskip

Applying Proposition 3.2 to receptive fields, N.~E.~Youngs 
in addition assumes that the relevant stimuli space $X$
forms a domain in $d$--dimensional real affine space, and that all non--empty intersections
of $U_i$ are convex. This assumption is also motivated by the orientation experiments,
mentioned in sec. 1 above.  She shows then that additional properties of $X$
can be extracted from $\C(\U)$, in particular, a lower bound for the embedding dimension $d$ and an information
about ``non--obvious'' relations between $U_i$ of the form
$$
\cap_{i\in\sigma}U_i\subset  \cup_{j\in\tau}U_j\ \roman{for\ certain}\ \sigma, \tau\  \roman{with}\  \sigma \cap \tau = \emptyset 
$$
([Yo14], p.~15). Algorithms for extracting this information are
given in terms of polynomial generators of the ideal in $\bold{F}_2[x_1,\dots ,x_n]$
defining the code.

\medskip

{\bf 4. Neural codes, semiotics, and modes of mathematical thinking.} A code,
as we defined it in the Introduction, is an analog of  {\it dictionary} in linguistics. It is
the list of minimal units of a language. Whereas theory of, say,
error--correcting codes does not invoke any a priori semantic notions,
semantics of neural codes was our main preoccupation here.
(Cf. similarly motivated discussion of genomics in [Ge93]
and Svante P\"a\"abo's diagnosis in his book ``Neanderthal Man: In Search of Lost Genomes'':
{\it ``The dirty little secret of genomics is that we still know next to nothing about how a genome translates into the particularities of a living and breathing individual.''})
\smallskip

Furthermore,  the main problem  addressed in  [CuIt08], [CuItVCYo13], [Yo14]
can be compared to the {\it decompilation} in computer science: translating a 
machine/assembly language code  back into a higher level language code.
In  semiotic parlance this means that spiking activity plays the role of {\it signifier}
whereas the nerve of the covering is the {\it signified} referring to a
 stimuli space, its covering, and individual positions in the stimuli space.

\smallskip

Below I  suggest that (a fragment of) the mathematical language
used in this description of the signifier of neural code  might have
a wider applicability. I will focus on the example of ``semantic space'',
 discussed in [Man08], sec.~2. In most variations, one starts with
 the assumption that  ``the space of meanings" (say, of words in a 
 language) is a certain set $X$.
 The next step consists in postulating additional structures
 (in  Bourbaki's sense) on it. 

\smallskip

(a) The most straightforward assumption is that meanings form
{\it subsets $U_i$} rather than points of $X$ so that we can imagine them as {\it a covering}
even before (or without) postulating a topology.

\smallskip

For example, P.~Guiraud ([Gui68]) suggests that  subsets of meaning $U_i$ have a natural
embedding into affine spaces $\bold{R}^d$, whose axes are marked by ``semes''.
Each seme axis corresponds to a  semantic  opposition (in the structuralist paradigm),
such as ``animate/inanimate'',
so that a pure "yes/no" picture would provide only an embedding into $\{-1,+1\}^d$ or else $\{0,1\}^d$.
However, a more realistic description of meanings would allow less localised positions
on the seme's axes. Say, on the axis ``animate/inanimate'', what place should be 
assigned to viruses?

\smallskip

(b) Discussing Zipf's Law, D.~Manin in [Man08] postulates furthermore that $X$ is endowed with 
{\it a measure}, and that $U_i$  are its measurable subsets. He uses this structure
in order to provide a mechanism generating Zipf's Law. 

\smallskip

On the other hand, in the generating model of Zipf's Law  developed in [Ma13] the key role
is played by the Kolmogorov complexity of the hypothetical {\it neural encoding}
of the semantic space rather than any specific structure of this space itself.
It is implicitly suggested that such  encodings are combinatorial objects,
exactly as in stimuli spaces encodings considered above. Therefore it would be extremely
interesting to study neural encodings of human languages from the perspective
of [CuIt08] and [CuItVCYo13].

\smallskip
(c) My last remark concerns a noteworthy parallel development in contemporary 
foundations of mathematics: Voevodsky's ``univalent foundations'', cf. [HTT13] and
a brief introduction in  [PeWa12]. Roughly speaking, in this model 
``continuous'' (homotopy type theory whose logical  description involves
simplicial sets, exactly as in the models of neural codes) belongs to the level of foundations, whereas ``discrete'' (set theory)
is a superstructure. During the whole XXth century, 
from the times of Georg Cantor and Felix Hausdorff, these two modes
of mathematical thinking were developing in the reverse order, from discrete to continuous:
cf. Figure 3 on p.~10 of [PeWa12].

\smallskip

This probably justifies the decompilation metaphor and the image of
``Foundations as Superstructure''  suggested in [Ma12].

\medskip

{\it Acknowledgements.} Preparing this survey, I benefited from e-mail exchanges 
with Nora Esther Youngs, Carina Curto,  Vladimir Itskov, Dmitri Manin, Ronald Brown.
\bigskip
\centerline{\bf References}

\medskip

[BrPo03] R.~Brown, T.~Porter. {\it Category theory and higher dimensional
algebra: potential descriptive tools in neuroscience. } Proceedings
of the International Conference on Theoretical Neurobiology, Delhi,
February 2003, ed.~by Nandini Singh, National Brain Research
Centre, Conference Proceedings 1 (2003) 80-92. arXiv:math/0306223

\medskip
[CuIt08] C.~Curto, V.~Itskov. {\it Cell Groups Reveal Structure of Stimulus Space.} 
PLoS Computational Biology, vol.~4, issue 10, October 2008, 13 pp. (available online).

\smallskip

[CuItVCYo13] C.~Curto, V.~Itskov, A.~Veliz-Cuba, N.~Youngs. {\it The neural ring:
An algebraic tool for analysing the intrinsic structure of neural codes.}
Bull. Math.~Biology, 75(9), pp. 1571--1611, 2013.
\smallskip

[DeSch01] E.~De Schutter. {\it Computational Neuroscience: More Math Is needed 
to Understand the Human Brain.} Mathematics Unlimited -- 2001 and Beyond.
Ed. by B\"jorn Engquist, Wilfried Scmidt. Springer 2001, pp. 381--391.

\smallskip

[Ge93] M.~S.~Gelfand. {\it Genetic Language -- Metaphor or Analogy?} Biosystems 30 (1--3), 1993,
pp. 277--288

\smallskip

[GeMa03] S.~I.~Gelfand, Yu.~I.~Manin. {\it Methods of homological algebra.} 2nd Edition,
Springer Verlag, 2003. xx+372 pp.

\smallskip

[GiIt14] Ch.~Giusti, V.~Itskov. {\it A No--Go Theorem for One--Layer Feedforward Networks.}
Neural Computation 26, 2014, pp. 2527--2540.

\smallskip

[Gui68] P.~Guiraud. {\it The semic matrices of meaning.} Social Science Information,
7(2), 1968, pp. 131--139.
\smallskip
[Ha02] A.~Hatcher. {\it Algebraic Topology.} CUP, Cambridge, 2002.

\smallskip

[HTT13] {\it Homotopy Type Theory: Univalent Foundations of Mathematics.} Institute for
Advanced Study, 2013.

\smallskip

[MaGaJoGoAshFraFri00] E.~A.~Maguire, D.~G.~Gadian, I.~S.~Jongsrude, C.~D.~Good,
J.~Ashburner, R.~S.~Frackowiak, C.~D.~Frith. {\it Navigation--related structural change
in the hippocampi of taxi drivers.} Proc.~Nat.~Ac.~Sci.~USA, vol.~97, no.~8, pp.~4398--4403,
2000.

\smallskip

[Man08]  D.~Yu.~Manin. {\it Zipf's Law and Avoidance of Excessive Synonymy.} 
Cognitive Science 32 (2008), 1075--1098. arXiv: 0710.0105 [cs.CL]
\smallskip

[Ma99] Yu.~I. Manin. {\it Classical computing, quantum computing,
and Shor's factoring algorithm.}  S\'eminaire Bourbaki, no. 862 (June 1999),
Ast\'erisque, vol 266, 2000, 375--404.
Preprint quant-ph/9903008.
\smallskip

[Ma11] Yu.~I. Manin. {\it A computability challenge: asymptotic bounds and isolated
error--correcting codes.} In: WTCS 2012 (Calude Festschrift), Ed. by
M.J. Dinneen et al.,
Lecture Notes in Computer Sci., 7160, pp. 174--182, 2012. Preprint
arXiv:1107.4246.
\smallskip
[Ma12] Yu.~I.~Manin. {\it Foundations as superstructure. (Reflections of a practicing
mathematician).} In:
``Philosophy, Mathematics, Linguistics: Aspects of Interaction. Proc.
of the
International Sci. Conference, St. Petersburg, Euler Int. Math.
Institute, May 22--25.
St. Petersburg 2012, pp. 98--111. Preprint arXiv:1205.6044

\smallskip

[Ma13] Yu.~I. Manin. {\it Zipf's law and L.~Levin's probability distributions.}  Functional Analysis and its Applications,
vol. 48, no. 2, 2014. DOI 10.107/s10688-014-0052-1.
Preprint arXiv:
1301.0427
\smallskip

[MaMar12] Yu.~I. Manin, M.~Marcolli.  {\it Kolmogorov complexity and the asymptotic bound for error-correcting
codes.}   Journ. of Diff.~Geom., 97, 2014, pp. 91--108.  Preprint arXiv:1203.0653
\smallskip

[PeWa12] \'A.~Pelayo, M.~A.~Warren. {\it Homotopy type theory and Voevodsky's univalent
foundations.} arXiv:1210.5658, 48 pp.

\smallskip
[SiGh07] V.~de Silva, R.~Ghrist. {\it Homological Sensor Networks.} Notices of the AMS,
vol.~54, no.~1, pp.~10--17, 2007.

\smallskip
[VlaNoTsfa07] S.~G.~Vladut, D.~Yu.~Nogin, M.~A.~Tsfasman. {\it Algebraic geometric codes: basic notions.} Mathematical Surveys and Monographs, 139. American Mathematical Society, Providence, RI, 2007.

\smallskip

[Ya00] V.~V.~Yashchenko (Ed.) et al. {\it Introduction to Cryptography.} (In Russian). MCCE, 2000.
\smallskip

[Yo14] N.~E.~Youngs. {\it The neural ring: using algebraic geometry to analyse neural rings.}
arXiv:1409.2544 [q-bio.NC], 108 pp.

\enddocument